\documentclass[12pt]{article}
\usepackage{amsmath}
\usepackage{amssymb}
\usepackage{vatola}

\def\q{\quad}
\def\qq{\qquad}
\def\mod#1{\ (\text{\rm mod}\ #1)}
\def\t{\hbox}
\def\qtq#1{\q\t{#1}\q}
\def\f{\frac}
\def\e{\equiv}

\def\b{\binom}
\def\o{\omega}

\def\sls#1#2{(\f{#1}{#2})}
 \def\ls#1#2{\big(\f{#1}{#2}\big)}
\def\Ls#1#2{\Big(\f{#1}{#2}\Big)}

\let \pro=\proclaim
\let \endpro=\endproclaim
\begin{document}
\leftline{preprint: November 20, 2013}\par\q\par\q
 \centerline {\bf
Cubic congruences and sums involving $\binom{3k}k$}
$$\q$$
\centerline{Zhi-Hong Sun} $$\q$$ \centerline{School of Mathematical
Sciences, Huaiyin Normal University,} \centerline{Huaian, Jiangsu
223001, P.R. China} \centerline{E-mail: zhihongsun@yahoo.com}
\centerline{Homepage: http://www.hytc.edu.cn/xsjl/szh}

\abstract{Let $p$ be a prime greater than $3$ and let $a$ be a
rational p-adic integer. In this paper we try to determine
$\sum_{k=1}^{[p/3]}\binom{3k}ka^k\pmod p$, and real the connection
between cubic congruences and the sum
$\sum_{k=1}^{[p/3]}\binom{3k}ka^k$, where $[x]$ is the greatest
integer not exceeding $x$.  Suppose that $a_1,a_2,a_3$ are rational
p-adic integers, $P=-2a_1^3+9a_1a_2-27a_3$, $Q=(a_1^2-3a_2)^3$ and
$PQ(P^2-Q)(P^2-3Q)(P^2-4Q)\not\equiv 0\pmod p$. In this paper we
show that the number of solutions of the congruence
$x^3+a_1x^2+a_2x+a_3\equiv 0\pmod p$ depends only on
$\sum_{k=1}^{[p/3]}\binom{3k}k(\frac{4Q-P^2}{27Q})^k\pmod p$.
 Let $q$ be a
prime of the form $3k+1$ and so $4q=L^2+27M^2$ with $L,M\in\Bbb Z$.
When $p\not=q$ and $p\nmid L$, we establish congruences for
$\sum_{k=1}^{[p/3]}\binom{3k}k(\frac{M^2}q)^k$ and
$\sum_{k=1}^{[p/3]}\binom{3k}k(\frac{L^2}{27q})^k$ modulo p. As a
consequence, we show that $x^3-qx-qM\equiv 0\pmod p$ has three
solutions if and only if $p$ is a cubic residue of $q$.

\par\q
 }
 \endabstract
 \footnotetext[1] {The author is
supported by the Natural Sciences Foundation of China (grant no.
11371163).}

\section*{1. Introduction}
\par\q\par Congruences involving binomial coefficients are
interesting, and they are concerned with Fermat quotients, Lucas
sequences, Legendre polynomials, binary quadratic forms and cubic
congruences. Let $\Bbb Z$ be the set of integers, and for a prime
$p$ let $\Bbb Z_p$ denote the set of  those rational numbers whose
denominator is not divisible by $p$.
   Let $p>5$ be a prime. In [12] Zhao, Pan and Sun proved
   that
$$\sum_{k=1}^{p-1}\b{3k}k2^k\e \f 65((-1)^{(p-1)/2}-1)\mod p.$$ In
[10] Z.W. Sun investigated $\sum_{k=0}^{p-1}\b{3k}ka^k\mod p$ for
$a\in\Bbb Z_p$. He gave explicit congruences for $a=-4,\f 16,\f
17,\f 18,\f 19,\f 1{13},\f 38,\f 4{27}$.
\par Suppose that $p>3$ is a prime and $k\in\{0,1,\ldots,p-1\}$. It
is easy to see that $p\mid \b{3k}k$ for $\f p3<k<\f p2$ and
$\f{2p}3<k<p$. Thus, for any $a\in\Bbb Z_p$,
$$\sum_{k=1}^{p-1}\b{3k}ka^k\e
\sum_{k=1}^{[p/3]}\b{3k}ka^k+\sum_{k=(p+1)/2}^{[2p/3]}\b{3k}ka^k
\mod p,$$ where $[x]$ is the greatest integer not exceeding $x$.  In
[9] the author investigated congruences for
$\sum_{k=0}^{[p/3]}\b{3k}{k}a^k$ modulo $p$. In this paper we reveal
the connection between cubic congruences and  the sum
$\sum_{k=1}^{[p/3]}\b{3k}{k}a^k$. Let $\sls mp$ be the Legendre
symbol. For a prime $p>3$ and $a_0,a_1,a_2,a_3\in\Bbb Z_p$ with
$a_0\not\e 0\mod p$ let $N_p(a_0x^3+a_1x^2+a_2x+a_3)$ denote the
number of solutions of the cubic congruence
$a_0x^3+a_1x^2+a_2x+a_3\e 0\mod p$. It is well known (see for
example [2,4,5] and [7, Lemma 2.3]) that
$$N_p(x^3+a_1x^2+a_2x+a_3)=1\iff \Ls Dp=-1,\tag 1.1$$
where $D$ is the discriminant of $x^3+a_1x^2+a_2x+a_3$ given by
$$D=a_1^2a_2^2-4a_2^3-4a_1^3a_3-27a_3^2+18a_1a_2a_3.\tag 1.2$$
\par In this paper we prove the following typical results:
\par Let $p>3$ be a prime and $a_1,a_2,a_3\in\Bbb Z_p$.
Suppose $P=-2a_1^3+9a_1a_2-27a_3$, $Q=(a_1^2-3a_2)^3$ and
$PQ(P^2-Q)(P^2-3Q)(P^2-4Q)\not\e 0\mod p$. Then
$N_p(x^3+a_1x^2+a_2x+a_3)=3$ if and only if
$\sum_{k=1}^{[p/3]}\b{3k}k\ls{4Q-P^2}{27Q}^k\e 0\mod p$, and
$N_p(x^3+a_1x^2+a_2x+a_3)=1$ if and only if
$$x\e \f
P{3(a_1^2-3a_2)}\sum_{k=0}^{[p/3]}\b{3k}k\Ls{4Q-P^2}{27Q}^k
-\f{a_1}3\mod p$$ is a solution of the congruence
$x^3+a_1x^2+a_2x+a_3\e 0\mod p$.
\par Let $p>3$ be a prime, $n\in\Bbb Z_p$  and
$3n+2\not\e 0\mod p$. Then
$$\sum_{k=1}^{[p/3]}\b{3k}k(n^2(n+1))^k\e
\f{3(n+1)}{2(3n+2)}\Big(\Ls{(1+n)(1-3n)}p-1\Big)\mod p.$$
\par Let
$p$ be a prime of the form $3k+1$, $m\in\Bbb Z_p$ and $m\not\e
1,-2,-\f 12\mod p$. Then
$$\align &\sum_{k=1}^{[p/3]}\b{3k}k\f 1{(-3(m-1)(m+2))^k}
\\&\e\cases 0\mod p&\t{if $\ls{m-1}{m+2}^{\f{p-1}3}\e 1\mod p$,}
\\-\f 3{2m+1}\big(m+1+\ls{m-1}{m+2}^{\f{p-1}3}\big)\mod p&\t{if
$\ls{m-1}{m+2}^{\f{p-1}3}\not\e 1\mod p$.}
\endcases\endalign$$

\par Let $q$ be a prime of the form $3k+1$ and so $4q=L^2+27M^2$ with
$L,M\in\Bbb Z$ and $L\e 1\mod 3$. Let $p$ be a prime with
$p\not=2,3,q$ and $p\nmid L$. In this paper we determine
$$\sum_{k=1}^{[p/3]}\b{3k}k\f{M^{2k}}{q^k}\qtq{and}
\sum_{k=1}^{[p/3]}\b{3k}k\f{L^{2k}}{(27q)^k}\mod p.$$ As a
consequence, we show that $x^3-qx-qM\e 0\mod p$ has three solutions
if and only if $p$ is a cubic residue of $q$.
\par Let  $\o=\f{-1+\sqrt{-3}}2$.  For a prime $p>3$ let $\ls{a+b\omega}p_3$ be the cubic
Jacobi symbol defined in [6]. For $c\in\Bbb Z_p$ with $c^2+3\not \e
0\mod p$ and $r\in\{0,1,2\}$ following [6] we define $c\in C_r(p)$
if and only if $\ls{c+1+2\omega}p_3=\omega^r$. In this paper we
prove that $c\in C_r(p)$ depends only on $\sum_{k=1}^{[p/3]}
\b{3k}k\ls 4{9(c^2+3)}^k\mod p$. See Theorem 2.5.
\par For $a,b,c\in\Bbb Z$ and a prime $p$, if there are integers $x$ and $y$
  such that $p=ax^2+bxy+cy^2$, throughout this paper we briefly write that $p=ax^2+bxy+cy^2$.
\section*{2. Main results}

\par For any numbers $P$ and $Q$, let $\{U_n(P,Q)\}$ be
the Lucas sequence given by
$$U_0(P,Q)=0,\ U_1(P,Q)=1,\ U_{n+1}(P,Q)=PU_n(P,Q)-QU_{n-1}(P,Q)
\ (n\ge 1).$$ It is well known (see [11]) that
$$\aligned &U_n(P,Q)\\&=\cases \f
1{\sqrt{P^2-4Q}}\Big\{\Big(\f{P+\sqrt{P^2-4Q}}2\Big)^n-
\Big(\f{P-\sqrt{P^2-4Q}}2\Big)^n\Big\} &\t{if $P^2-4Q\not=0$,}
\\n\ls P2^{n-1}&\t{if $P^2-4Q=0$.}
\endcases\endaligned\tag 2.1$$
Let $U_n=U_n(P,Q)$. Using (2.1) we see that (see [11]) for any
positive integer $n$,
$$U_{n+1}U_{n-1}-U_n^2=-Q^{n-1}\qtq{and}
U_{2n+1}=U_{n+1}^2-QU_n^2.\tag 2.2$$
\par From now on we use $\sls ap$ to denote the Legendre symbol.
\pro{Lemma 2.1 ([9, Lemma 3.3])} Let $p>3$ be a prime and
$P,Q\in\Bbb Z_p$ with $PQ\not\e 0\mod p$. Then

$$U_{2[\f p3]+1}(P,Q)\e\cases -Q^{1-\f{p-\sls p3}3}
U_{\f{p-\sls p3}3-1}(P,Q)\mod p&\t{if}\ \ls{P^2-4Q}p=1,
\\-Q^{-\f{p-\sls p3}3}U_{\f{p-\sls p3}3+1}(P,Q)\mod p&\t{if}
\ \ls{P^2-4Q}p=-1.\endcases$$\endpro

 \pro{Lemma 2.2}  Let $p>3$ be
a prime and $P,Q\in\Bbb Z_p$ with $PQ\not\e 0\mod p$. Then
$$U_{\f{p-\sls p3}3}(P,Q)\e 0\mod p\qtq{implies}
U_{2[\f p3]+1}(P,Q)\e (-Q)^{[\f p3]}\mod p.$$ Moreover, if
$P^2-3Q\not\e 0\mod p$, then
$$U_{2[\f p3]+1}(P,Q)\e (-Q)^{[\f p3]}\mod p
\qtq{implies}U_{\f{p-\sls p3}3}(P,Q)\e 0\mod p.$$
\endpro
Proof. Set $U_m=U_m(P,Q)$ and $n=(p-\sls p3)/3$. Then $2[\f
p3]+1=p-n$.  We first assume $U_n\e 0\mod p$. If $p\e 1\mod 3$, by
(2.2) we have $U_{n+1}(PU_n-U_{n+1})/Q-U_n^2=-Q^{n-1}$ and so
$U_{n+1}^2\e Q^n\mod p$. Hence $U_{2n+1}=U_{n+1}^2-QU_n^2\e
Q^n=(-Q)^n\mod p$. If $p\e 2\mod 3$, by (2.2) we have
$(PU_n-QU_{n-1})U_{n-1}-U_n^2=-Q^{n-1}$ and so $U_{n-1}^2\e
Q^{n-2}\mod p$. Thus, $U_{2[\f p3]+1}=U_{2n-1}=U_n^2-QU_{n-1}^2\e
-Q\cdot Q^{n-2}=(-Q)^{[\f p3]}\mod p$.
\par Now we assume that $P^2-3Q\not\e 0\mod p$ and
$U_{2[\f p3]+1}(P,Q)\e (-Q)^{[\f p3]}\mod p$. We claim that
$P^2-4Q\not\e 0\mod p$. If $P^2-4Q\e 0\mod p$, by (2.1) we have
$$\align U_{2[\f p3]+1}&\e
U_{2[\f p3]+1}(P,P^2/4) =\Big(2\big[\f p3\big]+1\Big)\Ls P2^{2[\f
p3]}
\\&\e\Big(2\big[\f
p3\big]+1\Big)Q^{[\f p3]}\e \f 13(-Q)^{[\f p3]}\not\e (-Q)^{[\f
p3]}\mod p,\endalign$$ which contradicts the assumption. Hence
$P^2-4Q\not\e 0\mod p$. Suppose $p\e 1\mod 3$. By Lemma 2.1 we have
$$\aligned &U_{n-1}\e -Q^{2n-1}\mod p\qtq{for}\ls{P^2-4Q}p=1,
\\&U_{n+1}\e -Q^{2n}\mod p\qtq{for}\ls{P^2-4Q}p=-1.\endaligned\tag 2.3$$
When $\sls{P^2-4Q}p=-1$ we have
$$ Q^n\e U_{2n+1}=U_{n+1}^2-QU_n^2
\e Q^{4n}-QU_n^2 \mod p.$$ As $Q^{3n}=Q^{p-1}\e 1\mod p$ we have
$Q^{4n}\e Q^n\mod p$. Thus $QU_n^2\e 0\mod p$ and so $U_n\e 0\mod
p$. When $\sls{P^2-4Q}p=1$ we have
$$\align Q^n&\e U_{2n+1}=U_{n+1}^2-QU_n^2=(PU_n-QU_{n-1})^2-QU_n^2
\\&\e (PU_n+Q^{2n})^2-QU_n^2
=U_n((P^2-Q)U_n+2PQ^{2n})+Q^{4n}\mod p.\endalign$$ As
 $Q^{4n}\e Q^n\mod p$ we have
$$U_n((P^2-Q)U_n+2PQ^{2n})\e 0\mod p.$$
If $P^2\e Q\mod p$, as $p\nmid PQ$ we have $U_n\e 0\mod p$. Now
assume $P^2-Q\not\e 0\mod p$.
 If $U_n\e -\f{2PQ^{2n}}{P^2-Q}\mod p$, then
$$U_{n+1}=PU_n-QU_{n-1}\e
-\f{2P^2Q^{2n}}{P^2-Q}+Q^{2n}=\f{Q+P^2}{Q-P^2}Q^{2n}\mod p.$$ Hence
$$\align -Q^{n-1}&=U_{n+1}U_{n-1}-U_n^2
\e \f{Q+P^2}{Q-P^2}Q^{2n}(-Q^{2n-1})-\f{4P^2Q^{4n}}{(P^2-Q)^2}
\\&=\f{Q^{4n-1}}{(P^2-Q)^2}(P^4-Q^2-4P^2Q)\mod p.\endalign$$
As $Q^{4n-1}\e Q^{n-1}\mod p$ we must have
$$P^4-Q^2-4P^2Q\e -(P^2-Q)^2\mod p.$$
That is, $2P^2(P^2-3Q)\e 0\mod p.$ This contradicts the assumption.
Thus, $(P^2-Q)U_n+2PQ^{2n}\not\e 0\mod p$ and so $U_n\e 0\mod p$.
\par
Now we assume $p\e 2\mod 3$. As $U_{2n-1}=U_{2[\f p3]+1}\e (-Q)^{[\f
p3]}=-Q^{n-1}\mod p$, by Lemma 2.1 we have
$$\aligned &U_{n-1}\e Q^{2n-2}\mod p\qtq{for}\ls{P^2-4Q}p=1,
\\&U_{n+1}\e Q^{2n-1}\mod p\qtq{for}\ls{P^2-4Q}p=-1.\endaligned
\tag 2.4$$ When $\sls{P^2-4Q}p=1$ we have
$$ -Q^{n-1}\e U_{2n-1}=U_n^2-QU_{n-1}^2
\e U_n^2-Q^{4n-3} \mod p.$$ As $Q^{4n-3-(n-1)}=Q^{3n-2}=Q^{p-1}\e
1\mod p$ we have $Q^{4n-3}\e Q^{n-1}\mod p$. Thus $U_n^2\e 0\mod p$
and so $U_n\e 0\mod p$. When $\sls{P^2-4Q}p=-1$ we have
$$\align -Q^{n-1}&\e U_{2n-1}=U_n^2-QU_{n-1}^2
=U_n^2-Q\Ls{PU_n-U_{n+1}}Q^2\e U_n^2-\f{(PU_n-Q^{2n-1})^2}Q
\\&=-\f {U_n((P^2-Q)U_n-2PQ^{2n-1})}Q-Q^{4n-3}\mod p.\endalign$$ As
 $Q^{4n-3}\e Q^{n-1}\mod p$ we have
$$U_n((P^2-Q)U_n-2PQ^{2n-1})\e 0\mod p.$$ If $P^2-Q\e 0\mod p$, as $p\nmid PQ$ we have
$U_n\e 0\mod p$. Now assume $P^2-Q\not\e 0\mod p$.
 If $U_n\e \f{2PQ^{2n-1}}{P^2-Q}\mod p$, then
$$ U_{n-1}=\f{PU_n-U_{n+1}}Q\e
\f PQ\cdot \f{2PQ^{2n-1}}{P^2-Q}-Q^{2n-2}
=Q^{2n-2}\f{2P^2-(P^2-Q)}{P^2-Q}\mod p.$$
 Hence
$$\align -Q^{n-1}&=U_{n+1}U_{n-1}-U_n^2
\e Q^{4n-3}\Big(\f{2P^2-(P^2-Q)}{P^2-Q}-\f{4P^2Q}{(P^2-Q)^2}\Big)
\\&\e \f{Q^{n-1}}{(P^2-Q)^2}(-(P^2-Q)^2+2P^2(P^2-3Q))\mod p.\endalign$$
This yields
 $2P^2(P^2-3Q)\e 0\mod p,$ which contradicts the assumption.
Thus, $(P^2-Q)U_n-2PQ^{2n-1}\not\e 0\mod p$ and so $U_n\e 0\mod p$.
\par Summarizing all the above we prove the lemma.

 \pro{Lemma 2.3 ([9,
(3.1)]} Let $p>3$ be a prime and $P,Q\in\Bbb Z_p$ with $PQ\not\e
0\mod p$. Then
$$U_{2[\f p3]+1}(P,Q)\e (-Q)^{[\f p3]}\sum_{k=0}^{[p/3]}
\b{3k}k\Ls{P^2}{27Q}^k\mod p.$$
\endpro
 \pro{Lemma 2.4}  Let $p>3$ be a prime and $P,Q\in\Bbb Z_p$ with
$PQ(P^2-3Q)(P^2-4Q)\not\e 0\mod p$. Then the following statements
are equivalent:
\par $(\t{\rm i})$ $U_{(p-\sls p3)/3}(P,Q)\e 0\mod p,$
\par $(\t{\rm ii})$
$\sum_{k=1}^{[p/3]}\b{3k}k\ls{P^2}{27Q}^k\e 0\mod p,$
\par $(\t{\rm iii})$ The congruence $x^3-3Qx-PQ\e 0\mod p$ has three
solutions.
\endpro
Proof. By Lemmas 2.2 and 2.3,
$$\align U_{\f{p-\sls p3}3}(P,Q)\e 0\mod p&\iff
U_{2[\f p3]+1}(P,Q)\e (-Q)^{[\f p3]}\mod p\\& \iff
\sum_{k=1}^{[p/3]}\b{3k}k\Ls{P^2}{27Q}^k\e 0\mod p.\endalign$$ Thus
(i) is equivalent to (ii). By [8, (7.4)] or [6, Corollary 6.3], (i)
is equivalent to (iii).

 \pro{Theorem 2.1} Let $p>3$ be
a prime and $a\in\Bbb Z_p$ with $a\not\e 0,\f 19,\f 1{27},\f
4{27}\mod p$. Then the following statements are equivalent:
\par $(\t{\rm i})$ $ \sum_{k=1}^{[p/3]}\b{3k}ka^k\e
0\mod p$,
\par $(\t{\rm ii})$ $U_{(p-\sls p3)/3}(9a,3a)\e 0\mod p,$
\par $(\t{\rm iii})$ $\ls{27a-2+3\sqrt{81a^2-12a}}2^{(p-\sls
p3)/3}\e 1\mod p$,
\par $(\t{\rm iv})$ $ax^3-x-1\e 0\mod p$ has three solutions,
\par $(\t{\rm v})$ $ \sum_{k=1}^{[p/3]}\b{3k}k\ls{4-27a}{27}^k\e
0\mod p$,
\par $(\t{\rm vi})$ $(27a-4)x^3+3x+1\e 0\mod p$ has three solutions.

\endpro Proof. Taking $P=9a$ and $Q=3a$ in
Lemma 2.4 we see that (i) and (ii) are equivalent. By (2.1),
$$U_{\f{p-\sls p3}3}(9a,3a)\e 0\mod p
\iff \Ls{9a+\sqrt{(9a)^2-4\cdot 3a}}{9a-\sqrt{(9a)^2-4\cdot 3a}}
^{\f{p-\sls p3}3}\e 1\mod p.$$ As
$$\f{9a+\sqrt{81a^2-12a}}{9a-\sqrt{81a^2-12a}}
=\f{(9a+\sqrt{81a^2-12a})^2}{12a}=\f {27a-2+3\sqrt{81a^2-12a}}2,$$
we see that (ii) is equivalent to (iii). For $x=3ay$ we see that
$$x^3-3\cdot 3ax-9a\cdot 3a=(3ay)^3-9a\cdot
3ay-27a^2=27a^2(ay^3-y-1).$$ Thus, $x^3-3\cdot 3ax-9a\cdot 3a\e
0\mod p$ has three solutions if and only if $ay^3-y-1\e 0\mod p$ has
three solutions. Hence applying Lemma 2.4 we see that (ii) is
equivalent to (iv). It is clear that
$$\align &\f{27(\f 4{27}-a)-2+3\sqrt{81(\f 4{27}-a)^2-12(\f
4{27}-a)}}2\cdot \f{27a-2+3\sqrt{81a^2-12a}}2
\\&=\f{2-27a+3\sqrt{81a^2-12a}}2\cdot \f{27a-2+3\sqrt{81a^2-12a}}2
=-1.\endalign$$ Thus,
$$\align &\Ls{27(\f 4{27}-a)-2+3\sqrt{81(\f 4{27}-a)^2-12(\f
4{27}-a)}}2^{\f{p-\sls p3}3}\e 1\mod p\\&\iff
\Ls{27a-2+3\sqrt{81a^2-12a}}2^{\f{p-\sls p3}3}\e 1\mod p.\endalign$$
Since $\f 4{27}-a\not\e 0,\f 19,\f 4{27}$ and (iii) is equivalent to
(i) and (iv), using the above we see that (iii) is equivalent to (v)
and that
$$\align &ax^3-x-1\q\t{has three solutions}
\\&\iff(\f 4{27}-a)x^3-x-1\e 0\mod p\q \t{has three solutions}
\\&\iff(\f 4{27}-a)(3x)^3-3x-1\e 0\mod p\q \t{has three solutions}
\\&\iff (27a-4)x^3+3x+1\e 0\mod p\q \t{has three solutions.}
\endalign$$
Thus (iv) and (vi) are equivalent. Now the proof is complete.

\pro{Lemma 2.5} Let $p>3$ be a prime and $a\in\Bbb Z_p$ with
$a\not\e 0,\f 19,\f 1{27},\f 4{27}\mod p$. Then the cubic congruence
$(27a-4)x^3+3x+1\e 0\mod p$ has one and only one solution if and
only if $x\e \sum_{k=0}^{[p/3]}\b{3k}ka^k\mod p$ is a solution of
the congruence.
\endpro
Proof. As $27a-4\not\e 0\mod p$ and
$$(27a-4)^2((27a-4)x^3+3x+1)=((27a-4)x)^3
+3(27a-4)\cdot (27a-4)x+(27a-4)^2,$$ we see that
$$N_p((27a-4)x^3+3x+1)=N_p(x^3+3(27a-4)x+(27a-4)^2).$$
By (1.2), the discriminant of $x^3+3(27a-4)x+(27a-4)^2$ is
$27^2a(4-27a)^3$. If $N_p((27a-4)x^3+3x+1)=1$, by (1.1) and the
above we must have $\sls{a(4-27a)}p=-1$. Now applying [9, Theorem
3.10] we see that the unique solution of $(27a-4)x^3+3x+1\e 0\mod p$
is given by $x\e \sum_{k=0}^{[p/3]}\b{3k}ka^k\mod p$. Conversely,
suppose that $x\e \sum_{k=0}^{[p/3]}\b{3k}ka^k\mod p$ is a solution
of the congruence $(27a-4)x^3+3x+1\e 0\mod p$. If
$N_p((27a-4)x^3+3x+1)=3$, by Theorem 2.1 we have
$\sum_{k=0}^{[p/3]}\b{3k}ka^k\e 1\mod p$. But $x\e 1\mod p$ is not a
solution of the congruence $(27a-4)x^3+3x+1\e 0\mod p$. This
contradicts the assumption. Hence  $N_p((27a-4)x^3+3x+1)=1$. This
proves the lemma.

\pro{Theorem 2.2} Let $p>3$ be a prime and $a_1,a_2,a_3\in\Bbb Z_p$.
Suppose $P=-2a_1^3+9a_1a_2-27a_3$, $Q=(a_1^2-3a_2)^3$ and
$PQ(P^2-Q)(P^2-3Q)(P^2-4Q)\not\e 0\mod p$. Then
$N_p(x^3+a_1x^2+a_2x+a_3)=3$ if and only if
$\sum_{k=1}^{[p/3]}\b{3k}k\ls{4Q-P^2}{27Q}^k\e 0\mod p$, and
$N_p(x^3+a_1x^2+a_2x+a_3)=1$ if and only if
$$x\e \f
P{3(a_1^2-3a_2)}\sum_{k=0}^{[p/3]}\b{3k}k\Ls{4Q-P^2}{27Q}^k
-\f{a_1}3\mod p$$ is a solution of the congruence
$x^3+a_1x^2+a_2x+a_3\e 0\mod p$.
\endpro
Proof. Set $x=\f P{3(a_1^2-3a_2)}y-\f{a_1}3$. Then
$x^3+a_1x^2+a_2x+a_3=-\f P{27}(-\f{P^2}Qy^3+3y+1).$ Thus,
$$\align &N_p(x^3+a_1x^2+a_2x+a_3)\\&=N_p\Big(-\f{P^2}Qx^3+3x+1\Big)
=N_p\Big(\Big(27\cdot\f{4Q-P^2}{27Q}-4\Big)x^3+3x+1\Big).\endalign$$
 Hence, applying Theorem 2.1 we see that
$$\align& x^3+a_1x^2+a_2x+a_3\e 0\mod p\q\t{has three solutions}
\\&\iff \Big(27\cdot\f{4Q-P^2}{27Q}-4\Big)x^3+3x+1\e
0\mod p\q\t{has three solutions}
\\&\iff \sum_{k=1}^{[p/3]}\b{3k}k\Ls{4Q-P^2}{27Q}^k\e 0\mod
p.\endalign$$
\par  From the above and Lemma 2.5 we see that
$$\align &N_p(x^3+a_1x^2+a_2x+a_3)=1
\\&\iff N_p\Big(-\f{P^2}Qy^3+3y+1\Big)=1
\\&\iff y\e
\sum_{k=0}^{[p/3]}\b{3k}k\Ls{4Q-P^2}{27Q}^k\mod p \ \t{satisfying} \
-\f{P^2}Qy^3+3y+1\e 0\mod p
\\&\iff x\e \f P{3(a_1^2-3a_2)}
\sum_{k=0}^{[p/3]}\b{3k}k\Ls{4Q-P^2}{27Q}^k-\f {a_1}3\mod p
\\&\q\qq\t{satisfying}\q x^3+a_1x^2+a_2x+a_3\e 0\mod p.
\endalign$$
  This completes the proof.

\pro{Theorem 2.3} Let $p>3$ be a prime, $n\in\Bbb Z_p$ and
$3n+2\not\e 0\mod p$. Then
$$\sum_{k=1}^{[p/3]}\b{3k}k(n^2(n+1))^k\e
\f{3(n+1)}{2(3n+2)}\Big(\Ls{(1+n)(1-3n)}p-1\Big)\mod p.$$
\endpro
Proof. Clearly the result is true for $n\e 0,-1\mod p$. Now assume
$n(n+1)\not\e 0\mod p$. If $3n-1\e 0\mod p$, by Lemma 2.3 we have
$$\align \sum_{k=1}^{[p/3]}\b{3k}k(n^2(n+1))^k&\e \sum_{k=1}^{[p/3]}\b{3k}k
\Ls 4{27}^k\e (-1)^{[\f p3]}U_{2[\f p3]+1}(2,1)-1 \\&=(-1)^{[\f
p3]}\Big(2\big[\f p3\big]+1\Big)-1\e -\f 23\e
-\f{3(n+1)}{2(3n+2)}\mod p.\endalign$$ Thus the result is true for
$n\e \f 13\mod p$. From now on we assume $3n-1\not\e 0\mod p$. It is
clear that $27n^2(n+1)-4=(3n+2)^2(3n-1)$. Set $a=n^2(n+1)$. Then
$a(27a-4)\not\e 0\mod p$. By [8, (7.4)],
$$\align &U_{\f{p-\sls p3}3}(9a,3a)\e 0\mod p
\\&\iff x^3-9ax-27a^2\e 0\mod p\q\t{has three solutions}
\\&\iff (3ax)^3-9a\cdot 3ax-27a^2\e 0\mod p\q\t{has three solutions}
\\&\iff ax^3-x-1\e 0\mod p\q\t{has three solutions}
\\&\iff (n+1)(nx)^3-nx-n\e 0\mod p\q\t{has three solutions}
\\&\iff (n+1)x^3-x-n\e 0\mod p\q\t{has three solutions}
\\&\iff (x-1)((n+1)(x^2+x)+n)\e 0\mod p\q\t{has three solutions}
\\&\iff (x-1)\Big(\big(x+\f 12\big)^2+\f{3n-1}{4(n+1)}\Big)
\e 0\mod p\q\t{has three solutions}
\\&\iff \Ls{(1+n)(1-3n)}p=1.
\endalign$$
Hence, if $\sls{(1+n)(1-3n)}p=1$, then $U_{(p-\sls p3)/3}(9a,3a)\e
0\mod p$ and so $U_{2[\f p3]+1}(9a,3a)\e 0\mod p$ by Lemma 2.2.
Applying Lemma 2.3 we find that $\sum_{k=1}^{[p/3]}\b{3k}ka^k\e
0\mod p$. \par Now we assume $\sls{(1+n)(1-3n)}p=-1$.
  As
$3n\not\e 0,-2\mod p$ we see that
$\sls{a(4-27a)}p=\sls{n^2(n+1)(3n+2)^2(1-3n)}p
=\sls{(1+n)(1-3n)}p=-1$. By [9, Theorem 3.10], $x\e
\sum_{k=0}^{[p/3]}\b{3k}k(n^2(n+1))^k\mod p$ is the unique solution
of the congruence $(3n+2)^2(3n-1)x^3+3x+1\e 0\mod p$. As
$$\align &(3n+2)^2(3n-1)x^3+3x+1
\\&=(3n-1)\Big(x+\f 1{3n+2}\Big)\Big\{\big((3n+2)x-\f
12\big)^2-\f{9(1+n)}{4(1-3n)}\Big\},
\endalign$$
we see that $x\e -\f 1{3n+2}\mod p$ is the unique solution of the
congruence $(3n+2)^2(3n-1)x^3+3x+1\e 0\mod p$. Hence
$\sum_{k=0}^{[p/3]}\b{3k}k(n^2(n+1))^k\e -\f 1{3n+2}\mod p$. This
completes the proof.

 \pro{Corollary 2.1} Let $p>3$ be a prime. Then
$$\align &
\sum_{k=1}^{[p/3]}\b{3k}k2^k\e \f 35\big((-1)^{\f{p-1}2}-1\big) \mod
p\qtq{for}p\not=5,
\\&\sum_{k=1}^{[p/3]}\b{3k}k(-4)^k\e \f 38\Big(\Ls {-7}p-1\Big)\mod p,
\\&\sum_{k=1}^{[p/3]}\b{3k}k12^k\e \f 9{16}\Big(\Ls
{-15}p-1\Big)\mod p,
\\&\sum_{k=1}^{[p/3]}\b{3k}k(-18)^k\e \f 37\Big(\Ls{-5}p-1\Big)
\mod p\qtq{for}p\not=7.\endalign$$
\endpro
Proof. Taking $n=1,-2,2,-3$ in Theorem 2.3 we deduce the result.

\pro{Corollary 2.2} Let $p>3$ be a prime. Then
$$\align &
\sum_{k=1}^{[p/3]}\b{3k}k36^k\e \f 6{11}\Big(\Ls{-2}p-1\Big) \mod
p\qtq{for}p\not=11,
\\&\sum_{k=1}^{[p/3]}\b{3k}k(-100)^k\e \f 6{13}\big((-1)^{\f{p-1}2}-1
\big)\mod p\qtq{for}p\not=13,
\\&\sum_{k=1}^{[p/3]}\b{3k}k\f 1{8^k}
\e \f 32\Big(\Ls 5p-1\Big)\mod p,
\\&\sum_{k=1}^{[p/3]}\b{3k}k\Ls 38^k
\e \f 9{14}\Big(\Ls {-3}p-1\Big)\mod p \qtq{for}p\not=7.
\endalign$$
\endpro
Proof. Taking $n=3,-5,-\f 12,\f 12$ in Theorem 2.3 we deduce the
result. \pro{Corollary 2.3} Let $p>3$ be a prime, $m\in\Bbb Z_p$ and
$m\not\e -3,9\mod p$. Then
$$\sum_{k=1}^{[p/3]}\b{3k}k\Ls{4(m-1)^2}{(m+3)^3}^k
\e\f 6{m-9}\Big(1-\Ls mp\Big)\mod p.$$
\endpro
Proof. Set $n=\f{1-m}{m+3}$. Then $n\not\e -1,-\f 23\mod p$ and
$m=\f{1-3n}{1+n}$. Now applying Theorem 2.3 we deduce the result.

\pro{Lemma 2.6 ([9, Theorem 3.3])} Let $p>3$ be a prime and
$a,b\in\Bbb Z_p$ with $ ab\not\e 0\mod p$. Then
$$\sum_{k=0}^{[p/3]}\b{3k}k\f{b^{2k}}{a^k}\e\cases
(-3a)^{[\f p3]+1}U_{\f{p-\sls p3}3-1}(9b,3a)\mod p&\t{\rm if}\
\sls{81b^2-12a}p=1,\\
-(-3a)^{[\f p3]}U_{\f{p-\sls p3}3+1}(9b,3a)\mod p&\t{\rm if}\
\sls{81b^2-12a}p=-1.\endcases$$\endpro \pro{Theorem 2.4} Let $p>3$
be a prime, $m\in\Bbb Z_p$ and $m\not\e 1,-2,-\f 12\mod p$.
\par $(\t{\rm i})$ If $p\e 1\mod 3$, then
$$\align &\sum_{k=1}^{[p/3]}\b{3k}k\f 1{(-3(m-1)(m+2))^k}
\\&\e\cases 0\mod p&\t{if $\ls{m-1}{m+2}^{\f{p-1}3}\e 1\mod p$,}
\\-\f 3{2m+1}\big(m+1+\ls{m-1}{m+2}^{\f{p-1}3}\big)\mod p&\t{if
$\ls{m-1}{m+2}^{\f{p-1}3}\not\e 1\mod p$.}
\endcases\endalign$$
\par $(\t{\rm ii})$ If $p\e 2\mod 3$, then
$$\align &\sum_{k=1}^{[p/3]}\b{3k}k\f 1{(-3(m-1)(m+2))^k}
\\&\e \f 1{2m+1}\Big\{(m-1)\Ls{m-1}{m+2}^{\f{p-2}3}+(m+2)
\Ls{m+2}{m-1}^{\f{p-2}3}\Big\}-1\mod p.\endalign$$
\endpro
Proof. As $81-12(-3(m-1)(m+2))=3^2(2m+1)^2$ and $[\f p3]+\f{p-\sls
p3}3=p-1-[\f p3]$, putting $a=-3(m-1)(m+2)$ and $b=1$ in Lemma 2.6
and then applying (2.1) we see that
$$\align &1+\sum_{k=1}^{[p/3]}\b{3k}k\f 1{(-3(m-1)(m+2))^k}
\\&\e (9(m-1)(m+2))^{[\f p3]+1}U_{\f{p-\sls p3}3-1}(9,-9(m-1)(m+2))
\\&=\f {(9(m-1)(m+2))^{[\f p3]+1}}{3(2m+1)}
\Big\{\Ls{9+3(2m+1)}2^{\f{p-\sls p3}3-1}-\Ls{9-3(2m+1)}2^{\f{p-\sls
p3}3-1}\Big\}
\\&\e \f{((m-1)(m+2))^{[\f p3]+1}}{2m+1}
\Big\{(m+2)^{\f{p-\sls p3}3-1}+(m-1)^{\f{p-\sls p3}3-1}\Big\}
\\&\e \f 1{2m+1}\Big\{(m-1)\Ls{m-1}{m+2}^{[\f
p3]}+(m+2)\Ls{m+2}{m-1}^{[\f p3]}\Big\} \mod p.
\endalign$$
If $p\e 1\mod 3$ and $\sls{m-1}{m+2}^{\f{p-1}3}\e 1\mod p$, from the
above we deduce that $$\sum_{k=1}^{[p/3]}\b{3k}k\f
1{(-3(m-1)(m+2))^k}\e \f 1{2m+1}(m-1+m+2)-1=0\mod p.$$ If $p\e 1\mod
3$ and $\sls{m-1}{m+2}^{\f{p-1}3}\not\e 1\mod p$, then
$1+\sls{m-1}{m+2}^{\f{p-1}3}+\sls{m-1}{m+2}^{-\f{p-1}3}\e 0\mod p$.
Thus, from the above we deduce that
$$\align&\sum_{k=1}^{[p/3]}\b{3k}k\f 1{(-3(m-1)(m+2))^k}
\\&\e \f 1{2m+1}
\Big(-(m+2)-3\Ls{m-1}{m+2}^{\f{p-1}3}\Big)-1 \\&=-\f
3{2m+1}\Big(m+1+\Ls{m-1}{m+2}^{\f{p-1}3}\Big)\mod p.\endalign$$ This
completes the proof.
 \pro{Corollary 2.4} Let $p>3$
be a prime, $m\in\Bbb Z_p$ and $(2m+1)^2\not\e 0,-3,9,-27\mod p$.
Then the congruence $x^3+3(m-1)(m+2)x+3(m-1)(m+2)\e 0\mod p$ has
three solutions if and only if $p\e 1\mod 3$ and
$\sls{m-1}{m+2}^{\f{p-1}3}\e 1\mod p$.
\endpro
Proof. Set $a=-\f 1{3(m-1)(m+2)}$. Then $a\not\e 0,\f 19,\f 1{27},\f
4{27}\mod p$. By Theorem 2.1,
$$\align &\sum_{k=1}^{[p/3]}\b{3k}k\f 1{(-3(m-1)(m+2))^k}
\e 0\mod p
\\&\iff \f 1{-3(m-1)(m+2)}x^3-x-1\e 0\mod p\q\t{has three solutions}
\\&\iff x^3+3(m-1)(m+2)x+3(m-1)(m+2)\e 0\mod p\q\t{has
three solutions}.\endalign$$ If $p\e 1\mod 3$, $t^2\e -3\mod p$
$(t\in\Bbb Z)$ and $\sls{m-1}{m+2}^{\f{p-1}3}\not\e 1\mod p$, then
$\sls{m-1}{m+2}^{\f{p-1}3}\e \f{-1\pm t}2\mod p$. As $(2m+1)^2\not\e
-3\mod p$ we have $2m+1\not\e \pm t\mod p$ and so $m+1\not\e \f{1\pm
t}2\mod p$. Thus $m+1+\sls{m-1}{m+2}^{\f{p-1}3}\not\e 0\mod p$.
Hence applying Theorem 2.4(i) we deduce the result.
\par Now we assume $p\e 2\mod 3$. As $\f{m+2}{m-1}\not\e 1\mod p$ we
have $\sls{m-1}{m+2}^{p-2}\not\e 1\mod p$ and so
$\sls{m-1}{m+2}^{\f{p-2}3}\not\e 1\mod p$. Since $\f{m-1}{m+2}\not\e
\pm 1\mod p$ we have $\sls{m-1}{m+2}^2\not\e 1\mod p$ and so
$\sls{m-1}{m+2}^{p+1}\not\e 1\mod p$. Hence $\sls
{m-1}{m+2}^{\f{p-2}3}\not\e \f{m+2}{m-1}\mod p$. Therefore $(m-1)
\sls {m-1}{m+2}^{\f{p-2}3}+(m+2)\sls {m+2}{m-1}^{\f{p-2}3}\not\e
2m+1\mod p$. This together with Theorem 2.4(ii) yields the result in
the case $p\e 2\mod 3$. The proof is now complete.

 \pro{Corollary 2.5} Let $p>3$ be a prime. Then
$$\sum_{k=1}^{[p/3]}\b{3k}k\f 1{6^k}
\e\cases 0\mod p&\t{if $3\mid p-1$ and $2^{\f{p-1}3}\e 1\mod p$,}
\\3\cdot 2^{\f{p-1}3}\mod p
&\t{if $3\mid p-1$ and $2^{\f{p-1}3}\not\e 1\mod p$,}
\\2^{\f{2p-1}3}
-2^{\f{p+1}3}-1\mod p&\t{if $p\e 2\mod 3$.}
\endcases$$
\endpro
Proof. Taking $m=-1$ in Theorem 2.4 we deduce the result.

\pro{Corollary 2.6} Let $p>5$ be a prime. Then
$$\align&\sum_{k=1}^{[p/3]}\b{3k}k\f 1{(-12)^k}
\\&\e\cases 0\mod p&\t{if $3\mid p-1$ and $2^{\f{p-1}3}\e 1\mod p$,}
\\-\f 35(2^{\f{p-1}3}+3)\mod p
&\t{if $3\mid p-1$ and $2^{\f{p-1}3}\not\e 1\mod p$,}
\\\f 15\cdot 2^{\f{p+1}3}(2^{\f{p+1}3}+1)-1\mod p&\t{if $p\e 2\mod 3$.}
\endcases\endalign$$
\endpro
Proof. Taking $m=2$ in Theorem 2.4 we deduce the result.

\pro{Corollary 2.7} Let $p>5$ be a prime with $p\not=13$. Then
$$\align&\sum_{k=1}^{[p/3]}\b{3k}k\f 1{(-120)^k}
\\&\e\cases 0\mod p&\t{if $3\mid p-1$ and $5^{\f{p-1}3}\e 1\mod p$,}
\\-\f 3{13}(5^{\f{p-1}3}+7)\mod p
&\t{if $3\mid p-1$ and $5^{\f{p-1}3}\not\e 1\mod p$,}
\\\f {10}{13}\cdot 5^{\f{p-2}3}(2\cdot 5^{\f{p-2}3}+1)-1\mod p&\t{if $p\e 2\mod 3$.}
\endcases\endalign$$
\endpro
Proof. Taking $m=6$ in Theorem 2.4 we deduce the result.

\pro{Corollary 2.8} Let $p>3$ be a prime, $c\in\Bbb Z_p$ and
$c\not\e 0,1,-1\mod p$. Then
$$\align&\sum_{k=1}^{[p/3]}\b{3k}k\Ls{(c+1)^2}{27c}^k
\\&\e \cases 0\mod p&\t{if $3\mid p-1$ and $c^{\f{p-1}3}\e 1\mod p$,}
\\\f 1{c-1}\big((c+1)c^{\f{p-1}3}-(c-2)\big)\mod p&\t{if $3\mid p-1$
and $c^{\f{p-1}3}\not\e 1\mod p$,}
\\\f 1{1-c}\cdot c^{\f{p+1}3}(1-c^{\f{p-2}3})-1\mod p&\t{if $p\e 2\mod
3$.}
\endcases\endalign$$\endpro
Proof. Taking $m=\f{1-2c}{1+c}$ in Theorem 2.4 we derive the result.
\par\q\par For two numbers $P$ and $Q$ let $\{V_n(P,Q)\}$ be defined by
$$V_0(P,Q)=2,\ V_1(P,Q)=P,\q V_{n+1}(P,Q)=PV_n(P,Q)-QV_{n-1}(P,Q)\
(n\ge 1).$$ It is well known that
$$V_n(P,Q)=\Ls{P+\sqrt{P^2-4Q}}2^n+\Ls{P-\sqrt{P^2-4Q}}2^n.$$
From [11] we know that
$$V_n(P,Q)=PU_n(P,Q)-2QU_{n-1}(P,Q)=2U_{n+1}(P,Q)-PU_n(P,Q).\tag 2.5$$

 \pro{Lemma 2.7 ([6, Corollary 6.1])} Let $p>3$
be a prime, $c\in\Bbb Z_p$ and $c(c^2+3)\not\e 0\mod p$. Then
$$U_{\f{p-\sls p3}3}(6,3(c^2+3))\e\cases 0\mod p&\t{if $c\in
C_0(p)$,}\\\f 1{2c}(-3(c^2+3))^{-[\f p3]}\mod p&\t{if $c\in
C_1(p)$,}
\\-\f 1{2c}(-3(c^2+3))^{-[\f p3]}\mod p&\t{if $c\in C_2(p)$}
\endcases$$ and
$$V_{\f{p-\sls p3}3}(6,3(c^2+3))\e\cases 2(3(c^2+3))^{-[\f p3]}
\mod p&\t{if $c\in C_0(p)$,}\\-(3(c^2+3))^{-[\f p3]}\mod p&\t{if
$c\in C_1(p)\cup C_2(p)$.}\endcases$$
\endpro

\pro{Theorem 2.5} Let $p>3$ be a prime, $c\in\Bbb Z_p$ and
$c(c^2+3)\not\e 0\mod p$. Then
$$ \sum_{k=1}^{[p/3]}\b{3k}k\Ls 4{9(c^2+3)}^k\e
\cases 0\mod p&\t{if $c\in C_0(p)$,}
\\-\f {3(c+1)}{2c}\mod p&\t{if $c\in C_1(p)$,}
\\-\f{3(c-1)}{2c}\mod p&\t{if $c\in C_2(p)$.}
\endcases$$\endpro
Proof. Let $a=c^2+3$ and $b=\f 23$. Then $81b^2-12a=-3\cdot 4c^2$.
By Lemma 2.6,
$$ \align&\sum_{k=0}^{[p/3]}\b{3k}k\Ls 4{9(c^2+3)}^k\\&\e
\cases -(3(c^2+3))^{\f{p-1}3+1}U_{\f{p-1}3-1}(6,3(c^2+3))\mod
p&\t{if $3\mid p-1$,}
\\(3(c^2+3))^{\f{p-2}3}U_{\f{p+1}3+1}(6,3(c^2+3))\mod
p&\t{if $3\mid p-2$.}\endcases\endalign$$ If $p\e 1\mod 3$, by (2.5)
 and Lemma 2.7 we have
 $$\align U_{\f{p-1}3-1}(6,3(c^2+3))&=\f
 1{6(c^2+3)}
 \Big(6U_{\f{p-1}3}(6,3(c^2+3))-V_{\f{p-1}3}(6,3(c^2+3))\Big)
\\&\e\cases -(3(c^2+3))^{-\f{p-1}3-1}\mod p&\t{if $c\in C_0(p)$,}
\\\f{c+3}{2c}(3(c^2+3))^{-\f{p-1}3-1}\mod p&\t{if $c\in C_1(p)$,}
\\\f{c-3}{2c}(3(c^2+3))^{-\f{p-1}3-1}\mod p&\t{if $c\in C_2(p)$.}
\endcases\endalign$$
 If $p\e 2\mod 3$, by (2.5)
 and Lemma 2.7 we have
 $$\align U_{\f{p+1}3+1}(6,3(c^2+3))&=3U_{\f{p+1}3}(6,3(c^2+3))
 +\f 12V_{\f{p+1}3}(6,3(c^2+3))
\\&\e\cases (3(c^2+3))^{-\f{p-2}3}\mod p&\t{if $c\in C_0(p)$,}
\\-\f{c+3}{2c}(3(c^2+3))^{-\f{p-2}3}\mod p&\t{if $c\in C_1(p)$,}
\\-\f{c-3}{2c}(3(c^2+3))^{-\f{p-2}3}\mod p&\t{if $c\in C_2(p)$.}
\endcases\endalign$$
Hence
$$\sum_{k=0}^{[p/3]}\b{3k}k\Ls 4{9(c^2+3)}^k
\e\cases 1\mod p&\t{if $c\in C_0(p)$,}
\\-\f{c+3}{2c}\mod p&\t{if $c\in C_1(p)$,}
\\-\f{c-3}{2c}\mod p&\t{if $c\in C_1(p)$.}
\endcases$$ This yields the result.
\pro{Corollary 2.9} Let $p>3$ be a prime, $c\in\Bbb Z_p$, $c\not\e
0,\pm 1\mod p$ and $c^2+3\not\e 0\mod p$. Then
$$c\in C_0(p)\iff \sum_{k=1}^{[p/3]}\b{3k}k\Ls 4{9(c^2+3)}^k\e 0\mod
p.$$\endpro \pro{Corollary 2.10} Let $p>3$ be a prime, $m\in\Bbb
Z_p$ and $(2m+1)(m^2+m+7)\not\e 0\mod p$. Then
$$ \sum_{k=1}^{[p/3]}\b{3k}k\f 1{(m^2+m+7)^k}\e
\cases  0\mod p&\t{if $\f{2m+1}3\in C_0(p)$,}
\\-\f{3(m+2)}{2m+1}\mod p&\t{if $\f{2m+1}3\in C_1(p)$,}
\\-\f{3(m-1)}{2m+1}\mod p&\t{if $\f{2m+1}3\in C_2(p)$.}
\endcases$$\endpro
 Proof. Set $c=\f{2m+1}3$. Then
$\f 4{9(c^2+3)}=\f 1{m^2+m+7}$. Now the result follows from Theorem
2.5. \pro{Corollary 2.11} Let $p>3$ be a prime. Then
$$\sum_{k=1}^{[p/3]}\b{3k}k\f
1{9^k}\e \cases 0\mod p&\t{if $p\e \pm 1,\pm 2\mod 9$,}
\\-3\mod p&\t{if $p\e \pm 4\mod 9$.}
\endcases$$
\endpro
Proof. As
$$\Ls{1+1+2\o}p_3=\Ls{-2\o^2}p_3=\Ls{\o}p_3^2=\Ls{\o}p_3^{-1}
=\o^{\f{\sls p3p-1}3},$$ we see that $1\in C_1(p)$ if and only if
$p\e \pm 4\mod 9$. Hence taking $m=1$ in Corollary 2.10 we obtain
the result.
 \pro{Corollary 2.12} Let $p>3$ be a prime.
\par $(\t{\rm i})$ If $p\e 1\mod 4$, then
$$\sum_{k=1}^{[p/3]}\b{3k}k\Ls 29^k\e 0\mod p\iff p=x^2+81y^2\ \t{or}
\ 2x^2+2xy+41y^2.$$
\par $(\t{\rm ii})$ If $p\e 1,3\mod{8}$, then
$$\sum_{k=1}^{[p/3]}\b{3k}k\Ls 49^k\e 0\mod p\iff p=x^2+162y^2\ \t{or}
\ 2x^2+81y^2.$$
\par $(\t{\rm iii})$ If $p\e 1,5,7,11\mod{24}$, then
$$ \sum_{k=1}^{[p/3]}\b{3k}k
\Big(-\f 4{27}\Big)^k\e 0\mod p\iff p=x^2+54y^2\ \t{or}\
2x^2+27y^2.$$
\par $(\t{\rm iv})$ If $p\e 1,2,4,8\mod{15}$, then
$$ \sum_{k=1}^{[p/3]}\b{3k}k
\Big(-\f 1{27}\Big)^k\e 0\mod p\iff p=x^2+135y^2\ \t{or}\
5x^2+27y^2.$$
\endpro
Proof. This is immediate from [6, Theorem 5.2] and Corollary 2.9.
\par For two integers $m$ and $n$ let $(m,n)$ be the greatest common
divisor of $m$ and $n$, and let $[m,n]$ be the least common multiple
of $m$ and $n$. Then we have:

\pro{Corollary 2.13} Let $p$ and $q$ be distinct primes greater than
$3$,
  $m,n\in\Bbb Z$, $(mn(m^2-n^2)(m^2+3n^2),pq)=1$ and $\sls p3p\e \sls q3q
  \mod {[9,m^2+3n^2]}$. Then
  $$\sum_{k=1}^{[p/3]}\b{3k}k\Ls {4n^2}{9(m^2+3n^2)}^k\e 0\mod
p\iff \sum_{k=1}^{[q/3]}\b{3k}k\Ls {4n^2}{9(m^2+3n^2)}^k\e 0\mod
q.$$
\endpro
Proof. Suppose $m+n(1+2\o)=\pm \o^r(1-\o)^s(a+b\o)$ with
$a,b,r,s\in\Bbb Z$ and $a+b\o\e 2\mod 3$. Then $a^2-ab+b^2\mid
m^2+3n^2$ and so $\sls p3p\e \sls q3q\mod {a^2-ab+b^2}$. By [8,
(1.1)-(1.2)],
$$\align \Ls{\f
mn+1+2\o}p_3&=\Ls{m+n(1+2\o)}p_3=\Ls{\o^r(1-\o)^s(a+b\o)}p_3
\\&=\Ls {\o}p_3^r\Ls{1-\o}p_3^s\Ls{a+b\o}p_3=\o^{\f{1-\sls p3p}3r}\cdot
\o^{\f{2(1-\sls p3p)}3s}\Ls{-\sls p3p}{a+b\o}_3
\\&=\o^{\f{1-\sls q3q}3r}\cdot
\o^{\f{2(1-\sls q3q)}3s}\Ls{-\sls q3q}{a+b\o}_3=\Ls{\f
mn+1+2\o}q_3.\endalign$$ Thus, $\f mn\in C_0(p)$ if and only if $\f
mn\in C_0(q)$. Now taking $c=\f mn$ in Corollary 2.9 and applying
the above we derive the result.

 \pro{Theorem 2.6} Let
$q$ be a prime of the form $3k+1$ and so $4q=L^2+27M^2$ with
$L,M\in\Bbb Z$ and $L\e 1\mod 3$. Let $p$ be a prime with
$p\not=2,3,q$ and $p\nmid L$. Then
$$\sum_{k=1}^{[p/3]}\b{3k}k\f{M^{2k}}{q^k}\e
\cases 0\mod p&\t{if $p^{\f{q-1}3}\e 1\mod q$,}
\\\f{-3-9M/L}2\mod p&\t{if $p^{\f{q-1}3}\e \f{-1+9M/L}2\mod q$,}
\\\f{-3+9M/L}2\mod p&\t{if $p^{\f{q-1}3}\e \f{-1-9M/L}2\mod q$.}
\endcases$$
\endpro
Proof. When $p\mid M$, by [6, Corollary 2.1] we have $p^{\f{q-1}3}
\e 1\mod q$. Thus the result is true. Now assume $p\nmid M$. Set
$c=\f L{3M}$. Then $c(c^2+3)\not\e 0\mod p$ and $\f
4{9(c^2+3)}=\f{M^2}q$. By Theorem 2.5,
$$\sum_{k=1}^{[p/3]}\b{3k}k\f{M^{2k}}{q^k}\e
\cases 0\mod p&\t{if $L/(3M)\in C_0(p)$,}
\\-\f 32(1+\f{3M}L)\mod p&\t{if $L/(3M)\in C_1(p)$,}
\\-\f 32(1-\f{3M}L)\mod p&\t{if $L/(3M)\in C_2(p)$.}
\endcases$$
From [6, Corollary 2.1] we know that for $i=0,1,2$,
$$p^{\f{q-1}3}\e \Ls{-1-L/(3M)}2^i\mod q
\iff \f L{3M}\in C_i(p).\tag 2.6$$
As $\f L{3M}\e -\f{9M}L\mod q$,
from the above we deduce the result.

 \pro{Corollary 2.14} Let $p>7$ be a prime.
Then
$$\sum_{k=1}^{[p/3]}\b{3k}k\f 1{7^k}\e
\cases 0\mod p&\t{if $p\e \pm 1\mod 7$,}
\\-6\mod p&\t{if $p\e \pm 2\mod 7$,}
\\3\mod p&\t{if $p\e \pm  4\mod 7$.}
\endcases$$
\endpro
Proof. As $4\cdot 7=1^2+27\cdot 1^2$, taking $q=7$ and $L=M=1$ in
Theorem 2.6 we deduce the result.
\par Similarly, from Theorem 2.6 we deduce the following results.
 \pro{Corollary 2.15} Let $p$ be a
prime with $p\not=2,3,5,13$. Then
$$\sum_{k=1}^{[p/3]}\b{3k}k\f 1{13^k}\e
\cases 0\mod p&\t{if $p\e \pm 1,\pm 5\mod {13}$,}
\\-\f{12}5\mod p&\t{if $p\e \pm 2,\pm 3\mod {13}$,}
\\-\f 35\mod p&\t{if $p\e \pm 4,\pm 6\mod {13}$.}
\endcases$$
\endpro
\pro{Corollary 2.16} Let $p$ be a prime with $p\not=2,3,7,19$. Then
$$\sum_{k=1}^{[p/3]}\b{3k}k\f 1{19^k}\e
\cases 0\mod p&\t{if $p\e \pm 1,\pm 7,\pm 8\mod {19}$,}
\\-\f 67\mod p&\t{if $p\e \pm 2,\pm 3,\pm 5\mod {19}$,}
\\-\f {15}7\mod p&\t{if $p\e \pm 4,\pm 6,\pm 9\mod {19}$.}
\endcases$$
\endpro

\pro{Corollary 2.17} Let $p$ be a prime with $p\not=2,3,11,37$. Then
$$\sum_{k=1}^{[p/3]}\b{3k}k\f 1{37^k}\e
\cases 0\mod p&\t{if $p\e \pm 1,\pm 6,\pm 8,\pm 10,\pm 11,\pm 14\mod
{37}$,}
\\-\f {12}{11}\mod p&\t{if $p\e \pm 2,\pm 9,\pm 12,\pm 15,
\pm 16,\pm 17\mod {37}$,}
\\-\f {21}{11}\mod p&\t{if $p\e \pm 3,\pm 4,\pm 5,\pm 7,\pm 13,
,\pm 18\mod {37}$.}
\endcases$$
\endpro
\pro{Theorem 2.7} Let $q$ be a prime of the form $3k+1$ and so
$4q=L^2+27M^2$ with $L,M\in\Bbb Z$ and $L\e 1\mod 3$. Let $p$ be a
prime with $p\not=2,3,q$ and $p\nmid L(L^2-9M^2)$. Then the
congruence $x^3-qx-qM\e 0\mod p$ has three solutions if and only if
$p$ is a cubic residue of $q$.
\endpro
Proof. When $p\mid M$, we have $L^2\e 4q\mod p$ and so $x^3-qx-qM\e
0\mod p$ has three solutions. On the other hand, by [6, Corollary
2.1] and Euler's criterion, $p$ is a cubic residue of $q$. Thus the
result is true in this case. Now we assume $p\nmid M$.
 By Theorems 2.1 and 2.6,
$$\align &\t{$x^3-qx-qM\e 0\mod p$ has three solutions}
\\&\iff \t{$ (Mx)^3-qMx-qM\e 0\mod p$ has three solutions}
\\&\iff\t{$\f {M^2}qx^3-x-1\e 0\mod p$ has three solutions}
\\&\iff \sum_{k=1}^{[p/3]}\b{3k}k\f{M^{2k}}{q^k}\e 0\mod p
\\&\iff p^{\f{q-1}3}\e 1\mod q
\\&\iff \t{$p$ is a cubic residues of $q$}.\endalign$$
\par As examples, if $p>3$ is a prime, then
$$\align &\t{$x^3-7x-7\e 0\mod p$ has three solutions}
\Leftrightarrow p=7\ \t{or}\ p\e \pm 1\mod 7,
\\&\t{$x^3-13x-13\e 0\mod p$ has three solutions}
\Leftrightarrow p=13\ \t{or}\ p\e \pm 1,\pm 5\mod {13},
\\&\t{$x^3-31x-62\e 0\mod p$ has three solutions}
\\&\qq\Leftrightarrow p=31\ \t{or $p$ is a cubic residue of $31$.}
\endalign$$

\pro{Theorem 2.8} Let $q$ be a prime of the form $3k+1$ and so
$4q=L^2+27M^2$ with $L,M\in\Bbb Z$ and $L\e 1\mod 3$. Let $p$ be a
prime with $p\not=2,3,q$ and $p\nmid M$. Then
$$\sum_{k=1}^{[p/3]}\b{3k}k\Ls{L^2}{27q}^k\e
\cases 0\mod p&\t{if $p^{\f{q-1}3}\e 1\mod q$,}
\\\f {-3\pm L/(3M)}2\mod p&\t{if $p^{\f{q-1}3}\e \f{-1\mp L/(3M)}2\mod q$.}
\endcases$$
\endpro
Proof. Set $c=-\f{9M}L$. Then $c(c^2+3)\not\e 0\mod p$ and $\f
4{9(c^2+3)}=\f{L^2}{27q}$. By [6, Proposition 2.2], $-\f{9M}L\in
C_i(p)$ if and only if $\f L{3M}\in C_i(p)$. Thus, from Theorem 2.5
we deduce that
$$\sum_{k=1}^{[p/3]}\b{3k}k\Ls{L^2}{27q}^k\e
\cases 0\mod p&\t{if $L/(3M)\in C_0(p)$,}
\\\f 12(3+\f L{3M})\mod p&\t{if $L/(3M)\in C_1(p)$,}
\\\f 12(3-\f L{3M})\mod p&\t{if $L/(3M)\in C_2(p)$.}
\endcases$$
This together with (2.6) yields the result.

\pro{Corollary 2.18 (Kummer, see [1, Theorem 10.10.5] or [3,
Corollaries 2.16 and 2.25])} Let $q$ be a prime of the form $3k+1$
and so $4q=L^2+27M^2$ with $L,M\in\Bbb Z$ and $L\e 1\mod 3$. Let $p$
be a prime with $p\not=2,3,q$ and $p\nmid M(L^2-81M^2)$. Then the
congruence $x^3-3qx-qL\e 0\mod p$ has three solutions if and only if
$p$ is a cubic residue of $q$.
\endpro
Proof. When $p\mid L$, from [6, Proposition 2.1 and Corollary 2.1]
we know that $0\in C_0(q)$ and so $p$ is a cubic residue of $q$.
Thus the result is true in this case. Now we assume that $p\nmid L$.
By Theorems 2.1 and 2.8,
$$\align &\t{$x^3-3qx-qL\e 0\mod p$ has three solutions}
\\&\iff \t{$ (\f L3x)^3-3q\cdot \f L3x-qL\e 0\mod p$ has three solutions}
\\&\iff\t{$\f {L^2}{27q}x^3-x-1\e 0\mod p$ has three solutions}
\\&\iff \sum_{k=1}^{[p/3]}\b{3k}k\Ls {L^2}{27q}^k\e 0\mod p
\\&\iff p^{\f{q-1}3}\e 1\mod q
\\&\iff \t{$p$ is a cubic residues of $q$}.\endalign$$
\par We remark that we prove Corollary 2.18 without  cyclotomic
numbers.
 \pro{Corollary 2.19}  Let $p$ be a prime with
$p\not=2,3,7$. Then
$$\sum_{k=1}^{[p/3]}\b{3k}k\f 1{189^k}\e
\cases 0\mod p&\t{if $p\e \pm 1\mod 7$,}
\\-\f 43\mod p&\t{if $p\e \pm 2\mod 7$,}
\\-\f 53\mod p&\t{if $p\e \pm 4\mod 7$.}
\endcases$$
\endpro
Proof. As $4\cdot 7=1^2+27\cdot 1^2$, taking $q=7$ and $L=M=1$ in
Theorem 2.8 we obtain the result.

\pro{Conjecture 2.1} Let $q$ be a prime of the form $3k+1$ and so
$4q=L^2+27M^2$ with $L,M\in\Bbb Z$ and $L\e 1\mod 3$. Let $p$ be a
prime with $p\not=2,3,q$ and $p\nmid LM$. Then
$$\sum_{\f p2<k<\f{2p}3}\b{3k}k\f{M^{2k}}{q^k}
\e\cases 0\mod p&\t{if $p^{\f{q-1}3}\e 1\mod q$,}
\\\pm\f{3M}L\mod p&\t{if $p^{\f{q-1}3}\e\f{-1\pm 9M/L}2\mod q.$}
\endcases$$
and
$$\sum_{\f p2<k<\f{2p}3}\b{3k}k\f{L^{2k}}{(27q)^k}
\e\cases 0\mod p&\t{if $p^{\f{q-1}3}\e 1\mod q$,}
\\\pm \f L{9M}\mod p&\t{if $p^{\f{q-1}3}\e\f{-1\pm L/(3M)}2\mod q.$}
\endcases$$
\endpro

\end{document}